\documentclass[conference]{IEEEtran}
\IEEEoverridecommandlockouts
\pdfoutput=1
\usepackage{cite}

\usepackage{amsmath,amssymb,amsfonts}
\usepackage{algorithm}
\usepackage{algpseudocode}
\usepackage{graphicx}
\usepackage{textcomp}
\usepackage{xcolor}

\def\BibTeX{{\rm B\kern-.05em{\sc i\kern-.025em b}\kern-.08em
    T\kern-.1667em\lower.7ex\hbox{E}\kern-.125emX}}
\begin{document}

\title{PID-inspired Continuous-time Distributed Optimization\\
}

\author{\IEEEauthorblockN{1\textsuperscript{st} Meng Tao}
\IEEEauthorblockA{\textit{School of Mathematics} \\
\textit{Southeast University}\\
Nanjing, China \\
mengtao@seu.edu.cn}
\and
\IEEEauthorblockN{2\textsuperscript{nd} Dongdong Yue}
\IEEEauthorblockA{\textit{School of Mathematics} \\
\textit{Southeast University}\\
Nanjing, China \\
yued@seu.edu.cn}
\and
\IEEEauthorblockN{3\textsuperscript{rd} Jinde Cao}
\IEEEauthorblockA{\textit{School of Mathematics} \\
\textit{Southeast University}\\
Nanjing, China \\
jdcao@seu.edu.cn}
}

\maketitle

\begin{abstract}
This paper proposes two novel distributed continuous-time algorithms inspired by PID control to solve distributed optimization problems. The algorithms are referred to as first-order and second-order, respectively, depend on the intrinsic dynamics of the agents in the network. Sufficient conditions are derived so that both algorithms converge exponentially over undirected connected graphs. Finally, numerical simulations illustrate the effectiveness and efficiency of the proposed algorithms.$ \footnote{The revised version of this work has been accepted by The 49th Annual Conference of the IEEE Industria1 Electronics Society.} $
\end{abstract}


\section{Introduction}
Distributed optimization (DO) aims to deal with the sum of local cost functions underlying a distributed interaction topology, which has drawn increasing attention in recent years \cite{nedic2018network,yang2019survey}. The DO algorithms, divided into continuous-time algorithms and discrete-time algorithms, have wide applications  including smart grid, machine learning and resource allocation (see \cite{boyd2011distributed} and the reference therein).

Continuous-time DO algorithms, inspired by system and control theory, exhibit some appealing features \cite{wang2011control,yang2016multi,zhang2014distributed}. For instance, the tedious work of characterizing sufficient small step sizes can be avoided in continuous-time algorithms, while this procedure is usually required in discrete-time algorithms. Proportional-Integral (PI) based multi-agent dynamics is one of the most important classes of continuous-time solvers for DO problems\cite{kia2015distributed,yi2018distributed,li2020input}. A significant feature of PI-type DO algorithms, in line with PI control theory, is the integral feedback action of the controlled error systems, i.e., the local consensus errors over the decision variables. Initially, the idea of PI has been formulated into the so-called saddle-point dynamics \cite{wang2011control,6578120}, where the  integral feedback variables are needed to be shared among the agents. The saddle-point dynamics has received a lot of attentions in recent years since it not only admits a PI structure, but also has well-established geometric interpretations \cite{yue2021distributed}. To save communication costs, a variant of saddle-point dynamics called modified-Lagrangian-based (MLB) has been proposed in \cite{kia2015distributed}, where the integral feedback variables are not needed to be shared but kept privacy. Note that zero-gradient-sum (ZGS) algorithm proposed in \cite{6129483} seems not belong to a PI-type algorithm at the first sight. However, interestingly, it has been shown recently in \cite{yu2021generalized} that MLB and ZGS can be recast into a generalized PI-type DO algorithm.\par
The application of the Proportional-Integral-Derivative (PID) control, as an extension of PI control, to distributed optimization offers significant advantages. PID control not only incorporates past information like PI control but also introduces future prediction capabilities. This enhanced performance has made PID control widely recognized and applied in various fields such as industrial engineering and multi-agent control\cite{borase2021review,ang2005pid}. However, there is a limited amount of research exploring the application of PID-type algorithms in distributed optimization scenarios. Motivated by this gap, we aim to leverage the PID concept in distributed optimization, harnessing its predictive capabilities and proven effectiveness, to improve the performance and adaptability of optimization processes. This approach holds great promise for addressing the challenges posed by distributed environments and achieving more efficient and effective optimization outcomes.

In fact, only one relevant literature \cite{zhu2022distributed} has made an early attempt, where the authors proposed a second-order PID-type algorithm. Nevertheless, we note that the algorithm proposed in [16] admits slow convergence rate due to the lack of the friction terms in the second-order dynamics (cf. \cite{yu2016gradient} for details). Clearly, the research on PID-type DO algorithms still awaits a breakthrough, which is crucial especially when PI-type DO algorithms have received more and more research attentions nowadays.\par 
In this paper, we proposed two algorithms, first-order and second-order respectively, for solving unconstrained DO problems over a connected graph. It is rigorously proved that both algorithms converge exponentially. Simulation results show that the proposed PID-type DO algorithms outperform existing PI-type algorithms\cite{zhang2014distributed,kia2015distributed,yi2018distributed,li2020input}, and the existing PID-type algorithm \cite{zhu2022distributed}, in terms of quicker response and faster convergence. \par
The remaining part of this article is organized as follows. Section II introduces the preliminaries. Two types of PID-based DO algorithms and convergence analysis will be presented in the section III . Section  IV presents two numerical simulations. The conclusion will be given in Section V.
\section{PRELIMINARY}

\subsection{Notations}
Let $\mathbf{0}, \mathbf{1} \in \mathbb{R}^{N} $ represent all-zero column vector and all-one column vector, respectively, and let $I_{N}$ represent the $N\times N$ identity matrix. The symbols $\otimes$  and $\|\cdot\|$ represent the Kronecker product vectors or matrices norm, respectively. Set col$\{x,y,\lambda\}\triangleq(x^{\mathrm{T}},y^{\mathrm{T}},\lambda^{\mathrm{T}})^{\mathrm{T}}$. Throughout the paper, a matrix $A\succ0$ represents $A$ is positive-definite, and $\rho(A)$ represents the spectral radius, $\lambda_{\max} (A)(\lambda_{\min} (A))$ is the largest(minimum) eigenvalue of matrix $A$. Let $\mathcal{G}=(\mathcal{V}, \mathcal{E})$ denotes an undirected graph with vertices set $\mathcal{V}=\{1 ,2, \ldots, n\}$ and edge set $\mathcal{E} \in \mathcal{V} \times \mathcal{V}$. The Laplacian matrix $L$ of a connected graph is symmetric positive semidefinite and satisfies $L \mathbf{1}=\mathbf{0} $ and $ \mathbf{1}^{\mathrm{T}} L=\mathbf{0}^{\mathrm{T}}$. 
\newtheorem{definition}{Definition}
\newtheorem{rem}{\bf Remark}
\begin{definition}\label{def1}\cite{6129483}
    A continuous differentiable $f(x):\mathbb{R}^{n}\to\mathbb{R}$ is $m$-strongly convex if, for any $x,y\in \mathbb{R}^{n}$, one of the following equivalent conditions hold:
    \begin{itemize}
            \item[I:]$\langle\nabla f(x)-\nabla f(y),x-y\rangle\geq m\|x-y\|^{2}$
            \item[II:]$\|\nabla f(x)-\nabla f(y)\|\geq m\|x-y\| $
    \end{itemize}
\end{definition}
\begin{definition}\label{def2}\cite{6129483,attouch2000heavy}
    A continuous differentiable $f(x):\mathbb{R}^{n}\to\mathbb{R}$ is $l$-smooth if, for any $x,y\in \mathbb{R}^{n}$, one of the following equivalent conditions hold:
    \begin{itemize}
            \item[I:]$\langle\nabla f(x)-\nabla f(y),x-y\rangle\leq l\|x-y\|^{2}$
            \item[II:]$\|\nabla f(x)-\nabla f(y)\|\leq l\|x-y\|$
            \item[III:]$\left\langle\nabla f(x)-\nabla f(y), x-y\right\rangle \geq\large{\frac{1}{l}}\left\|\nabla f(x)-\nabla f(y)\right\|^{2}$ 
    \end{itemize}
\end{definition}
\newtheorem{lem}{\bf Lemma}
\subsection{Problem Statement:}
\renewcommand{\IEEEQED}{\IEEEQEDopen}
Suppose there are $N$ agents interacting over an  undirected connected graph $\mathcal{G}$. Consider the unconstrained DO problem
 \begin{equation}\label{eq1}
     \min\limits_{z} f(z)=\textstyle\sum_{i=1}^{N}f_{i}(z),
 \end{equation}
 where $z\in \mathbb{R}^{n}$ is the decision variable, the global function $f$, which is the sum of private local objective function $f_{i}:\mathbb{R}^{n}\to\mathbb{R}$,
 satisfies the following assumption. 
 \newtheorem{ass}{Assumption}
\begin{ass}\label{ass1}
$f$ is continuous differentiable, $m$-strongly convex and $l$-smooth. 
\end{ass}
\begin{rem}
    The $m$-strong convexity assumption is weaker than the assumption that each local function is strongly convex \cite{6129483,li2020input,yu2021generalized}, i.e., the local cost function $f_{i}(x)$ may not be strongly convex and may even be nonconvex.
\end{rem}

\begin{lem}\label{lem0}\cite{yi2017formation}
    For an undirected and connected graph $\mathcal{G}$ with the associated Laplacian matrix $L$, there exists a matrix $\Gamma \in \mathbb{R}^{N \times N} \succ 0$  such that  $L \Gamma=\Gamma L=\Pi$, where  $\Pi=I_{N}-\frac{1}{N} \mathbf{1}_{N} \mathbf{1}_{N}^{\mathrm{T}}$. Moreover, $\Pi L=L\Pi=L$ and $\rho(\Pi)=1$.
\end{lem}\par
Under Assumption \ref{ass1}, the necessary and sufficient condition for optimality of problem (\ref{eq1}) is $\textstyle\sum_{i=1}^{N}\nabla f_{i}(z)=0$.\par
Denote $x_{i}\in \mathbb{R}^{n}$ as agent $i$'s local estimation of the global optimal solution and let $x=$col$(x_{1},\dots,x_{N})$, problem (\ref{eq1}) is equavalent to
\begin{align}\label{eq3}
    \min\limits_{x}f(x)=\textstyle\sum_{i=1}^{N}f_{i}(x_{i}),
    \text{    s.t. }x_{i}=x_{j}\quad\forall i,j,
\end{align}
where the constraints guarantee the consistency of the optimal solution. Then the optimal solution of problem (\ref{eq3}) satisfy
\begin{equation}\label{eq4}
    \textstyle\sum_{i=1}^{N}\nabla f_{i}(x_{i})=0,\quad x_{i}=x_{j}\quad\forall i,j.
\end{equation}\par

\section{MAIN RESULTS}

This section considers two distributed algorithms based on the idea of PID, and the convergences are analyzed separately.

\subsection{PID-Based First-order Distributed Optimization}
To solve problem (\ref{eq1}), consider the dynamics (\ref{eqq1}) for agent $i\in\mathcal{V}$.
\begin{algorithm}
        \caption{PID-Based First-order DO}
        \begin{itemize}
            \item \textbf{Initialization:} 
                \begin{itemize}
                    \item[1)] Choose constants $c_{k},k=1,2,3,4$.
                    \item[2)] Choose any $x_{i}(0)\in \mathbb{R}^{n}$, and $\lambda_{i}(0)\in \mathbb{R}^{n}$ such that $\sum_{i=1}^{N} \lambda_{i}(0)=0.$
                \end{itemize}
            \item \textbf{Dynamics for agent $i$:} 
        \end{itemize}
        \begin{subequations}\label{eqq1}
            \begin{align}
                &\dot x_{i} =-c_{1}\nabla f_{i}(x_{i})-\phi_{i},\label{eq5a}\\
                &\phi_{i}=c_{2}\mu_{i}+\lambda_{i}+c_{3}y_{i},\label{eq5b}\\
                &\mu_{i}=\textstyle\sum_{i=1}^{N}a_{ij}(x_{i}-x_{j}),\label{eq5c}\\
                &\dot \lambda_{i}=c_{4}\mu_{i},\label{eq5d}\\
                &y_{i}=\textstyle\sum_{i=1}^{N}a_{ij}(-c_{1}\nabla f_{i}(x_{i})-\phi_{i}+c_{1}\nabla f_{j}(x_{j})+\phi_{j}),\label{eq5e}
            \end{align}
        \end{subequations}
        \label{Alg1}
\end{algorithm}
Each agent in (\ref{eq5a}) initially estimates the optimal value using local gradient descent. However, due to variations in the functions $f_i$, a PID auxiliary term $\phi_i$ is introduced to account for the disparity in local gradients. This compensation mechanism helps the agents reach a consensus on the final optimal value. The equations (\ref{eq5c}-\ref{eq5e}) represent the P, I, D terms, respectively, with respect to the diffusive coupling input $\mu_{i}$. The scalars $c_{k}$, $k=1,2,3,4$ are positive parameters. To guarantee the optimality of the equilibrium point, it is necessary to satisfy the initial condition $\textstyle\sum_{i=1}^{N} \lambda_{i}(x) = 0$. A straightforward choice for the initial condition is to set $\lambda_{i}(0) = 0$ for all agents $i$. Algorithm \ref{Alg1} is a distributed algorithm since each agent solely relies on information from its neighboring agents.\par      
Denote $\lambda=$col$(\lambda_{1},\dots,\lambda_{N})$. By eliminating $\mu_{i}$ and $y_{i}$, the compact form of dynamic (\ref{eqq1}) is written as
\begin{subequations}\label{eq7}
    \begin{align}
        &\dot x=-c_{1}\nabla f(x)-c_{2}\mathcal{L}x-\lambda-c_{3}\mathcal{L}\dot x,\label{eq7a}\\
        &\dot \lambda=c_{4}\mathcal{L}x,\label{eq7b}
    \end{align}
\end{subequations}
where $\mathcal{L}\triangleq L\otimes I_{n}$.\par
Before proceeding further, note that the diffusive couplings in (\ref{eq5b}) bring algebraic loops. It is important to verify the well-posedness of the feedback interconnection.\par
Equation (\ref{eq7a}) can be rewritten as
$(I+c_{3}\mathcal{L})\dot x=-c_{2}\mathcal{L}x-c_{1}\nabla f(x)-\lambda$.
Notice that $(I+c_{3}\mathcal{L})$ should be nonsingular such that system (\ref{eq7}) can be rewritten in the following explicit form, ensuring the well-posedness of the feedback interconnection:
\begin{subequations}\label{eq9}
    \begin{align}
        &\dot x=(I+c_{3}\mathcal{L})^{-1}(-c_{1}\nabla f(x)-c_{2}\mathcal{L}x-\lambda),\label{eq9a}\\
        &\dot \lambda=c_{4}\mathcal{L}x.\label{eq9b}
    \end{align}
\end{subequations}
\begin{rem}
    Note from (\ref{eq7}) that there exists derivative feedback, which is similar to the algorithm proposed in \cite{li2020input}. The difference is that the algorithm \ref{Alg1} feeds back the gradient of the state, while the algorithm in \cite{li2020input} feeds back the gradient of the auxiliary variable (see Eq. (26) of \cite{li2020input}). Another observation from (6) is that the proposed algorithm has the same structure as the algorithm proposed in \cite{yu2021generalized}. The difference is that the 
    algorithm in \cite{yu2021generalized} requires a diagonal matrix in place of the matrix $(I+c_4L)$, which is clearly not satisfied.
\end{rem}
\begin{lem}\label{lem2}
    Under Assumption \ref{ass1}, if there exists an equilibrium point $(x^{*},\lambda^{*})$ of system (\ref{eq9}) that satisfies $\textstyle\sum_{i=1}^{N} \lambda_{i}^{*}=0$, where $x^{*}=$col$(x^{*}_{1},\dots,x^{*}_{N})$ and $\lambda^{*}=$col$(\lambda_{1}^{*},\dots,\lambda_{N}^{*})$, then $(x^{*},\lambda^{*})$ is unique and $x^{*}=\mathbf{1}_{N}\otimes z^{*}$, with $z^{*}$ being the optimal solution to problem (\ref{eq1}).
\end{lem}
\begin{IEEEproof}
See Appendix A.
\end{IEEEproof}
Next, the convergence result of the algorithm \ref{Alg1} is given.
\newtheorem{thm}{\bf Theorem}
\begin{thm}\label{th1}
    Under Assumption \ref{ass1}, with the initial condition $\textstyle\sum_{i=1}^{N} \lambda_{i}(x)=0$ and for any given $c_{i}>0$, $i=1,2,3,4$, the state $x$ in Algorithm \ref{Alg1} converges exponentially to the optimal point $x^{*}$.
\end{thm}
\begin{IEEEproof}
See Appendix B.
\end{IEEEproof}
\begin{rem}
    There are several ways to implement Algorithm \ref{Alg1} in a distributed form. For instance, adding a memory module in the feedback loop may eliminate algebraic loops \cite{abdelwahed2013reduced,gonzalez2011effect}. Although algebraic loops can be addressed using numerical methods, their introduction can significantly increase the complexity and pose challenges to the numerical stability of the system's solution process. This, in turn, may impact the convergence behavior and speed of the system. However, in the following, we will present a second-order algorithm based on the PID concept that effectively circumvents the need for algebraic loops.
\end{rem}

\subsection{PID-Based Second-order Distributed Optimization}

We propose the following second-order PID algorithm:
\begin{algorithm}
        \caption{PID-Based Second-order DO}
        \begin{itemize}
            \item \textbf{Initialization:} 
                \begin{itemize}
                    \item[1)] Choose constants $c_{k}$ satisfy the condition (\ref{eq23}).
                    \item[2)] Choose $x_{i}(0),v_{i}(0)\in \mathbb{R}^{n}$, and $\lambda_{i}(0)\in \mathbb{R}^{n}$ such that $\sum_{i=1}^{N} \lambda_{i}(0)=0.$
                \end{itemize}
            \item \textbf{Dynamics for agent $i$:} 
        \end{itemize}
        \begin{subequations}
            \begin{align}
                &\dot x_{i} =v_{i},\label{eq19a}\\
                &\dot v_{i}=-c_{1}\nabla f_{i}(x_{i})-c_{2}\mu_{i}-c_{3}\lambda_{i}-c_{4}w_{i}-c_{5}v_{i},\label{eq19b}\\                &\dot\lambda_{i}=\mu_{i},\label{eq19c}\\
                &\mu_{i}=\textstyle\sum_{i=1}^{N}a_{ij}(x_{i}-x_{j}),\label{eq19d}\\
                &w_{i}=\textstyle\textstyle\sum_{i=1}^{N}a_{ij}(v_{i}-v_{j}).\label{eq19e}
            \end{align}
        \end{subequations}
        \label{Alg2}
\end{algorithm}\par
By eliminating $w_{i}$ and $\mu_{i}$, the compact form of algorithm \ref{Alg1} can be rewritten as
\begin{subequations}\label{eq20}
    \begin{align}
        &\dot x=v,\label{eq20a}\\
        &\dot v=-c_{1}\nabla f(x)-c_{2}\mathcal{L}x-c_{3}\lambda-c_{4}\mathcal{L}v-c_{5}v,\label{eq20b}\\
        &\dot\lambda=\mathcal{L}x.\label{eq20c}
    \end{align}
\end{subequations}
Similarly, the three terms $c_{2}\mathcal{L}x$, $c_{3}\lambda$ and $c_{4}\mathcal{L}v$ of (\ref{eq20b}) correspond to the P-I-D terms. Algorithm \ref{Alg2} is a distributed algorithm since each agent solely relies on information from its neighboring agents. Besides, there are no algebraic loops in (\ref{eq20}).\par
\begin{rem}
    If we let $c_{4}=1$ and $c_{2}=c_{3}=c_{5}$, by introducing $\hat{v}=v+\mathcal{L}x$, the Algorithm \ref{Alg2} degenerates to the algorithm proposed in \cite{yu2016gradient}. In reference \cite{zhu2022distributed}, a novel second-order distributed optimization algorithm based on the PID approach is introduced. However, this algorithm differs from the one presented in this paper as it lacks the crucial friction term in the second-order dynamics. This absence of the friction term can lead to slower response times and even hinder convergence in \cite{zhu2022distributed}. Moreover, Note that compared to \cite{zhu2022distributed}, the PID-based DO algorithm proposed in this paper reduces the sharing of the dual variable $\lambda_{i}$ between agents. Therefore, the proposed algorithm can reduce communication overload.
\end{rem}
\begin{lem}\label{lema}
        Under Assumption \ref{ass1}, if there exists an equilibrium point $(x^{*},v^{*},\lambda^{*})$ of system (\ref{eq20}) that satisfies $\textstyle\sum_{i=1}^{N} \lambda_{i}^{*}=0$, where $x^{*}=$col$(x^{*}_{1},\dots,x^{*}_{N})$, $\lambda^{*}=$col$(\lambda_{1}^{*},\dots,\lambda_{N}^{*})$, and $v^{*}=$col$(v^{*}_{1},\dots,v^{*}_{N})$, then $(x^{*},v^{*},\lambda^{*})$ is unique and $x^{*}=\mathbf{1}_{N}\otimes z^{*}$, with $z^{*}$ being the optimal solution to problem (\ref{eq1}). Moreover, the following equation holds:
        \begin{subequations}\label{eq100}
            \begin{align}
                &\dot x^{*}=v^{*}\equiv 0,\label{eq100a}\\
                &\dot v^{*}=-c_{1}\nabla f(x^{*})-c_{3}\lambda^{*}\equiv 0,\label{eq10a}\\
                &\dot \lambda^{*}=\mathcal{L}x^{*}\equiv 0.\label{eq10b}
            \end{align}
        \end{subequations}
    \end{lem}
    \begin{IEEEproof}
        The analysis of Lemma 3 is similar to that of Lemma 1 and is omitted here.
    \end{IEEEproof}  
    
Next, we will give the convergence result of Algorithm \ref{Alg2}.
\begin{thm}\label{th2}
    Suppose $f$ is strictly convex and $l$-smooth, with the condition 
    \begin{equation}\label{eq23}
        (c_{1}l+c_{3}^{2}+c_{5}^{2}+\sqrt{\left(1+c_{2}^{2}+2c_{4}^{2}\right)\lambda_{\max}\left(L^{\mathrm{T}}L\right)})^{\frac{1}{2}}<\frac{\eta}{\gamma},
    \end{equation}
    every individual solution $x_{i}(t)$ exponentially converges to the unique global minimizer $x^{*}$, where $\eta$ and $\gamma$ are defined in (\ref{eq101}).
\end{thm}
\begin{IEEEproof}
See Appendix C.
\end{IEEEproof}
\begin{rem}
    The strict convexity of $f$ guarantees the existence
and uniqueness of optimal solution. Moreover, it is a weaker condition compared to strong convexity. It is worth mentioning that, differently
from \cite{yue2021distributed}, no additional assumptions on the form of $\nabla f_{i}$ are assumed here.
\end{rem}

The variation of Algorithm \ref{Alg2} can be conducted as the following Algorithm
\begin{align}\label{eq105}
    &\dot x=v,\nonumber\\
    &\dot v=-c_{1}\nabla f(x)-c_{2}\mathcal{L}x-c_{3}\mathcal{L}\lambda-c_{4}\mathcal{L}v-c_{5}v,\nonumber\\
    &\dot\lambda=\mathcal{L}x.
\end{align}

Next, we give the similar result to Theorem 2 by the following corollary.
\newtheorem{cor}{\bf Corollary}
\begin{cor}\label{cor1}
    Suppose $f$ is strictly convex and $l$-smooth, with the condition 
    \begin{equation}
        (c_{1}l+c_{5}^{2}+\sqrt{\left(1+c_{2}^{2}+c_{3}^{2}+2c_{4}^{2}\right)\lambda_{\max}\left(L^{\mathrm{T}}L\right)})^{\frac{1}{2}}<\frac{\eta}{\gamma},
    \end{equation}
    every individual solution $x_{i}(t)$ exponentially converges to the unique global minimizer $x^{*}$.
\end{cor}

The proof of Corollary \ref{cor1} is similar to Theorem \ref{th2}, and we omit here due to space limitations. Moreover, if we let $c_{5}=0$, by introducing $\Tilde{\lambda}=c_{2}^{2}c_{1}\lambda$, Algorithm (\ref{eq105}) degenerates to the algorithm proposed in \cite{zhu2022distributed}.
\begin{rem}
    In practical applications, the decision variables are often constrained, and the algorithm can be extended to distributed optimization problems with convex set constraints through the prox method.
\end{rem}
\section{NUMERICAL EXAMPLES}
This section provides two numerical examples to
verify the efficiency of the proposed algorithm.
\subsection{Example 1}

\begin{figure}[h]
  \begin{minipage}[t]{0.241\textwidth}
    \centering
    \includegraphics[width=\textwidth]{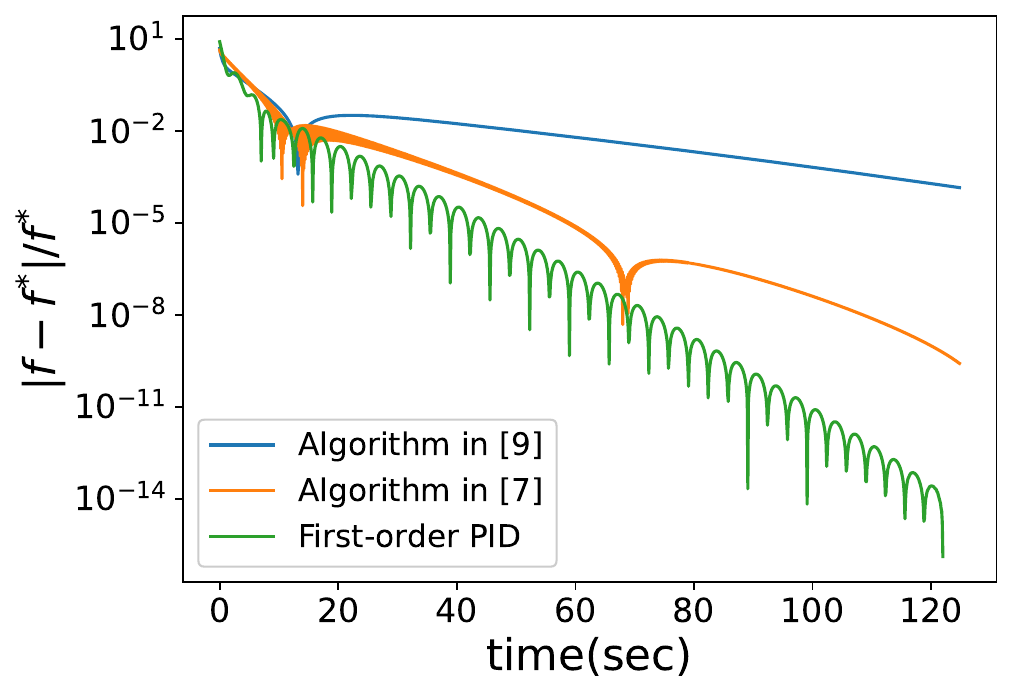}
    \caption{The relative error comparison of first-order algorithm (convex $f_{i}$)}
    \label{fig:example1}
  \end{minipage}
  \hfill
  \begin{minipage}[t]{0.241\textwidth}
    \centering
    \includegraphics[width=\textwidth]{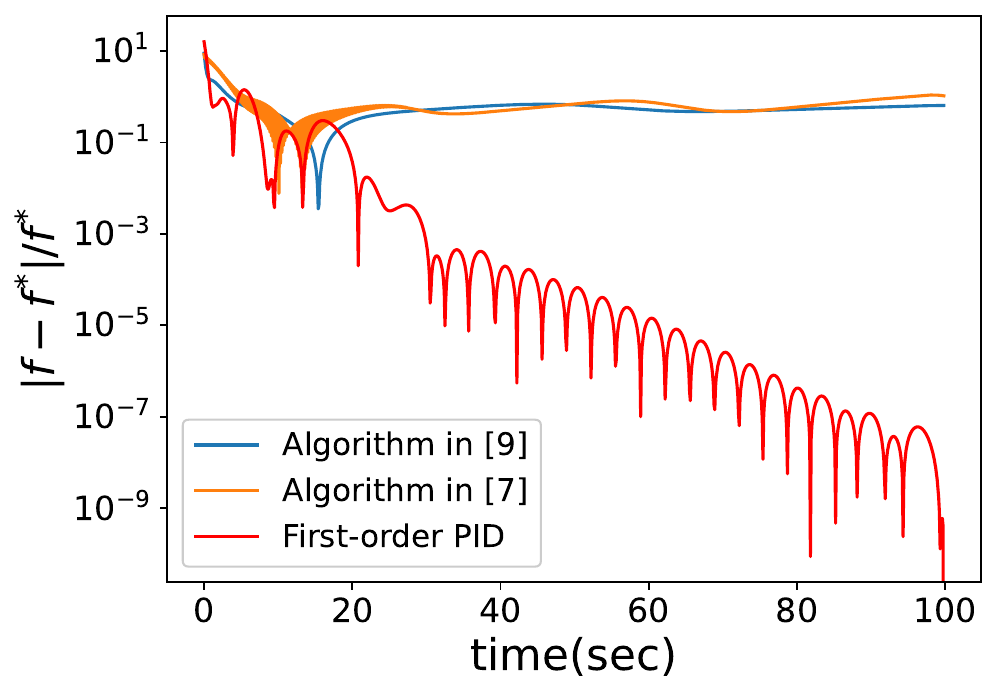}
    \caption{The relative error comparison of first-order algorithm (nonconvex $f_{i}$)}
    \label{fig:example2}
  \end{minipage}
\end{figure}

Consider the quadratic programming DO problem: $\min\limits_{x} \textstyle\sum_{i=1}^{N}\frac{1}{2} x^{\mathrm{T}} Q_{i} x+q_{i}^{\mathrm{T}} x$, where $Q_{i} \in \mathbb{R}^{n \times n}$ is a symmetric positive semi-definite matrix, $q_{i} \in \mathbb{R}^{n}$. Let $Q_{i}$ be generated by the uniform distribution of $[0,1]$, $q_{i}$ be generated by the uniform distribution of $[-5,5]$, taking $N = 4$ and $n = 10$. Under the ring topology of the network, choose $c_{1}=0.8,c_{2}=2.9,c_{3}=5,c_{4}=5$, the numerical results compared with the algorithm proposed in \cite{li2020input} and the algorithm proposed in \cite{kia2015distributed}, are shown in the Fig. \ref{fig:example1}, which validates the faster convergence of Algorithm \ref{Alg1}. Furthermore, in order to verify the effectiveness of the algorithm when the local objective function is non-convex, redefine $f_{1}(z)=\frac{1}{2} x^{\mathrm{T}} Q_{1} x+q_{1}x+\sin(z),f_{2}(z)=\frac{1}{2} x^{\mathrm{T}} Q_{2} x+q_{2}x-\sin(z), f_{3}(z)=\frac{1}{2} x^{\mathrm{T}} Q_{3} x+q_{3}x-5\cos(z), f_{4}=\frac{1}{2} x^{\mathrm{T}} Q_{4} x+q_{4}x+5\cos(z)$, the convergence result of the first-order PID algorithm is shown in Figure 2.
\subsection{Example 2}
Similar to the convex case of example 1, consider the quadratic programming problem above to verify the effectiveness of the proposed second-order dynamical system. Under the ring graph with 20 agents, the state $x_{i}\in\mathbb{R}^{7}$, choose $c_{1}=0.14,c_{2}=0.65,c_{3}=0.156,c_{4}=0.52,c_{5}=0.52$. The experimental results compared with the algorithm proposed in \cite{zhu2022distributed}, the algorithm proposed in \cite{yi2018distributed} and the algorithm proposed in \cite{zhang2014distributed} are shown in the Fig. \ref{fig:example2}. It can be seen that the proposed algorithm \ref{Alg2} outperforms the rest algorithms.
\begin{figure}[h]
    \centering
    \includegraphics[width=0.395\textwidth]{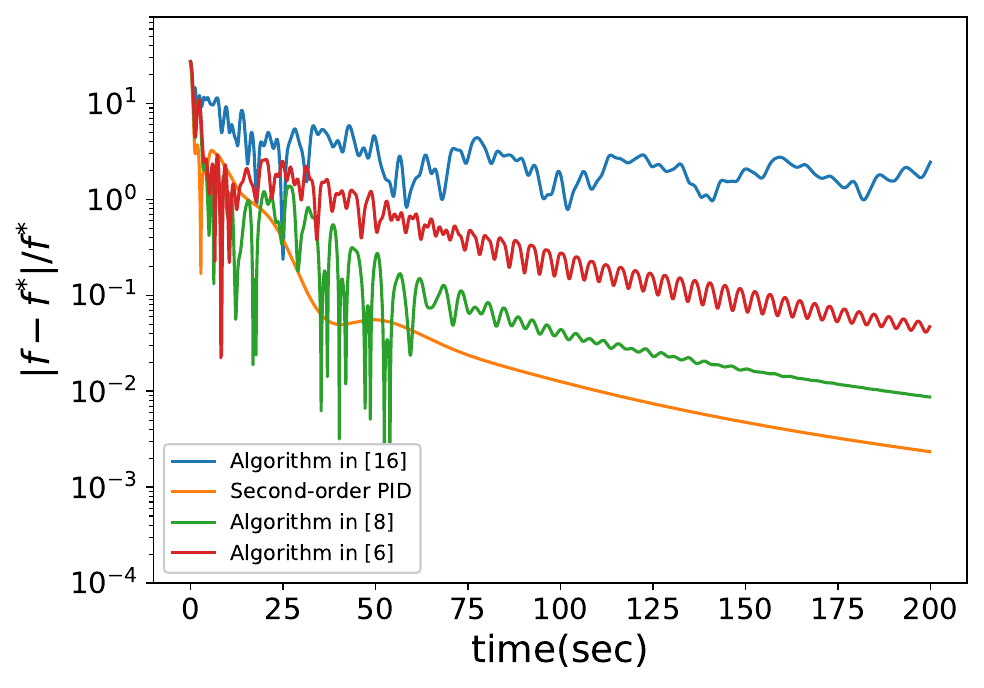}
    \caption{The relative error comparison of second-order algorithm}
    \label{fig:example2}
\end{figure}\par

\section{CONCLUSION}
This paper successfully apply the ideas of Proportional-Integral-Derivative (PID) control to solve distributed optimization (DO) problems. Two continuous-time algorithms, first-order and second-order respectively, are proposed and analyzed. Both algorithms achieve exponentially convergence in an connected graph. Finally, the effectiveness of the algorithm was verified by numerical examples. For future work, one may consider DO problems with directed topologies.

\section{APPENDIX}

\subsection{Appendix A}
Suppose $(x^{*},\lambda^{*})$ is the equilibrium point which satisfies \begin{subequations}\label{eq10}
        \begin{align}
            &\dot x^{*}=(I+c_{3}\mathcal{L})^{-1}(-c_{1}\nabla f(x^{*})-\lambda^{*})\equiv 0,\label{eq10a}\\
            &\dot \lambda^{*}=c_{4}\mathcal{L}x^{*}\equiv 0,\label{eq10b}
        \end{align}
\end{subequations}
where the term $-c_{2}\mathcal{L}x^{*}$ in (\ref{eq10a}) is zero and omitted since (\ref{eq10b}) implied $\mathcal{L}x^{*}\equiv 0$, which results in $x_{i}^{*}=x_{j}^{*}$ $\forall i,j$. Note that (\ref{eq10a}) can implies 
\begin{equation}\label{eq11}
    \dot x^{*}=-c_{1}\nabla f(x^{*})-\lambda^{*}\equiv 0.
\end{equation}
Next, multiplying (\ref{eq11}) by $(\mathbf{1}_{N}\otimes I_{n})^{\mathrm{T}}$ from the left, note that $(\mathbf{1}_{N} \otimes I_{m})^{\mathrm{T}} \lambda^{*}=\textstyle\sum_{i \in \mathcal{N}} \lambda_{i}^{*}=0$, one has
\begin{equation*}
    \begin{aligned}
        & -(\mathbf{1}_{N} \otimes I_{m})^{\mathrm{T}} c_{1} \nabla f(x^{*})-(\mathbf{1}_{N} \otimes I_{m})^{\mathrm{T}} \lambda^{*} \nonumber\\
        = & -c_{1} \textstyle\sum_{i \in \mathcal{N}} \nabla f_{i}(x_{i}^{*}) \equiv \mathbf{0},\nonumber
    \end{aligned}
\end{equation*}
which satisfies the optimality condition (\ref{eq4}).

\subsection{Appendix B}
Based on Lemma \ref{lem2}, we introduce the auxiliary state 
\begin{equation}\label{eq12}
    \theta \triangleq g(x)+w[\lambda+c_{1} \nabla {f}(x^{*})],
\end{equation}
where $w>0$, $g(x)\triangleq(I+c_{3}\mathcal{L})(x-x^{*})$. Consequently, the closed-loop system in (\ref{eq9}) can be rewritten as
\begin{align}\label{eq13}
    &\dot{x}=(I+c_{3}\mathcal{L})^{-1}(-c_{1} h(x)-c_{2} \mathcal{L} x-\frac{1}{w} \theta+\frac{1}{w} g(x)), \nonumber\\
    &\dot{\theta}=-c_{1} h(x)-c_{2} \mathcal{L} x-\frac{1}{w} \theta+\frac{1}{w} g(x)+c_{4} w \mathcal{L} x,
\end{align}
where $h(x)\triangleq\nabla f(x)-\nabla f(x^{*})$. Due to the smoothness of $f$,$$\|h(x)\|\leq l\|x-x^{*}\|
\text{, and }\|g(x)\|\leq \beta\|x-x^{*}\|,$$ with $\beta=c_{3}\lambda_{max}(L)+1$. We construct the following Lyapunov function:
\begin{equation}\label{eq14}
    V(x,\theta)=\frac{1}{2}(x-x^{*})^{\mathrm{T}}(c_{3}\mathcal{L}+I)(x-x^{*})+\frac{q}{2}\theta^{\mathrm{T}}\Gamma_{n}\theta,
\end{equation}
where $\Gamma_{n}\triangleq\Gamma\otimes I_{n}$, and $q>0$ is to be given. Since $(c_{3}\mathcal{L}+I)$ and $\Gamma$ are positive definite matrix,
\begin{align}
    &V(x, \theta) \geq \frac{1}{2}\lambda_{\min }(c_{3}\mathcal{L}+I)\|x-x^{*}\|^{2}+\frac{q}{2} \lambda_{\min }(\Gamma)\|\theta\|^{2},\nonumber\\
    &V(x, \theta) \leq \frac{1}{2}\lambda_{\max}(c_{3}\mathcal{L}+I)\|x-x^{*}\|^{2}+\frac{q}{2} \lambda_{\max }(\Gamma)\|\theta\|^{2}.\nonumber
\end{align}
Since $\sum_{i=1}^{N} \lambda_{i}(0)=0$ and $\frac{\mathrm{d}}{\mathrm{d} t}(\mathbf{1}_{N}^{\mathrm{T}} \otimes I_{n}) \lambda(t)=0$, we have $(\mathbf{1}_{N}^{\mathrm{T}} \otimes I_{n}) \lambda(t)=(\mathbf{1}_{N}^{\mathrm{T}} \otimes I_{n})(\lambda-\lambda^{*})=0$.
By defining $\Gamma_{n}$, one can obtain that
    \begin{equation}\label{eq15}
        \theta^{\mathrm{T}} \Gamma_{n} \mathcal{L} x=\theta^{\mathrm{T}}(x-x^{*})-\frac{1}{N} g^{\mathrm{T}}(\mathbf{1}_{N} \mathbf{1}_{N}^{\mathrm{T}} \otimes I_{n})(x-x^{*}).
    \end{equation}
    For simplicity, we will drop the augument $x$ in $h(x)$ and $g(x)$. Then, the derivative of $V$ along the solution of (\ref{eq13}) satisfies
    \begin{align*}
        \dot{V}= & (x-x^{*})^{\mathrm{T}}(c_{3}\mathcal{L}+I)\dot x+q \theta^{\mathrm{T}} \Gamma_{n} \dot{\theta} \\
        =&(x-x^{*})^{\mathrm{T}}(-c_{1} h-c_{2} \mathcal{L} x-\frac{1}{w} \theta+\frac{1}{w} g)\\&+q \theta^{\mathrm{T}} \Gamma_{n}\left[-c_{1} h-\frac{1}{w} (\theta- g)\right]\\&+q \theta^{\mathrm{T}} \Gamma_{n}(c_{4} w-c_{2}) \mathcal{L}(x-x^{*}) \\
        \leq & -c_{2} x^{\mathrm{T}} \mathcal{L} x-c_{1} m\|x-x^{*}\|^{2}-\frac{1}{w} \theta^{\mathrm{T}}(x-x^{*}) \\
        & +\frac{1}{w} g^{\mathrm{T}}(x-x^{*})+q \theta^{\mathrm{T}} \Gamma_{n}\left[-c_{1} h-\frac{1}{w}(\theta-g)\right] \\
        & +q \theta^{\mathrm{T}} \Gamma_{n}(c_{4} w-c_{2}) \mathcal{L}(x-x^{*}).
    \end{align*}
    By selecting $w$, $q>0$ such that $c_{4} w-c_{2}>0$ and $q=\frac{1}{w(c_{4} w-c_{2})}$ and applying (\ref{eq15}), we have
    \begin{align}\label{eq16}
        \dot{V} \leq & -c_{2} x^{\mathrm{T}} \mathcal{L} x-\left(c_{1} m-\frac{2 \beta}{w}\right)\|x-x^{*}\|^{2} \nonumber\\
        & -\frac{q}{w} \lambda_{\min }(\Gamma)\|\theta\|^{2}-q \theta^{\mathrm{T}} \Gamma_{n}\left(c_{1} h-\frac{g}{w}\right).
    \end{align}
    According to Young's inequality, there exist positive scalars $s_1$ and $s_2$, such that
    \begin{align*}
        \dot{V} \leq & -\left(c_{1} m-\frac{2 \beta}{w}\right)\|x-x^{*}\|^{2}-\frac{q}{w} \lambda_{\min }(\Gamma)\|\theta\|^{2} \\
        & +\frac{q \lambda_{\max }(\Gamma)}{w}\left[s_{1} l+s_{2} \beta\right]\|\theta\|^{2} \\
        & +\lambda_{\max }(\Gamma) q\left[\frac{w l c_{1}^{2}}{4 s_{1}}+\frac{\beta}{4 s_{2} w}\right]\|x-x^{*}\|^{2}.
    \end{align*}
    Meanwhile, for any given $c_{k}>0, k=1,2,3,4$, there exist sufficiently small $s_{j}, j=1,2$ and sufficiently large $w$, such that
    \begin{equation}\label{eq17}
        \{\begin{array}{l}
            \frac{2 \beta}{w}+\lambda_{\max }(\Gamma)\left(\frac{c_{1}^{2} l}{4 s_{1}\left(c_{4} w-c_{2}\right)}+\frac{\beta}{4 s_{2} w^{2}\left(c_{4} w-c_{2}\right)}\right) \leq \frac{c_{1} \alpha_{f}}{2}, \\
            \lambda_{\max }(\Gamma)\left(s_{1} l+s_{2} \beta\right) \leq \frac{\lambda_{\min }(\Gamma)}{2}.
        \end{array}
    \end{equation}
    Consequently, from $q=\frac{1}{w(c_{4} w-c_{2})}$, (\ref{eq17}) implies that
    \begin{equation}\label{eq18}
        \dot{V}\leq-\frac{c_{1}l}{2}\,\|x-x^{*}\|^{2}-\frac{q\lambda_{\mathrm{min}}(\Gamma)}{2w}\,\|\theta\|^{2}\leq-\sigma_{V}V,
    \end{equation}
    with $\sigma_{V}\triangleq\mathrm{min}\{\frac{c_{1}l}{\beta},\frac{\lambda_{\mathrm{min}}(\Gamma)}{w\lambda_{\mathrm{max}}(\Gamma)}\}$. Therefore, $x$ converges to $x^{*}$ exponentially.

    \subsection{Appendix C}
    Let ${\bar{x}}=x-x^{*},{\bar{v}}=v-v^{*},{\bar{\lambda}}=\lambda-\lambda^{*}$, $\xi=(\bar{x}^{\mathrm{T}},\bar{v}^{\mathrm{T}},\bar{\lambda}^{\mathrm{T}})^{\mathrm{T}}$, and define
\begin{equation*}
		A=\left(
		\begin{matrix}
			c_{4}\mathcal{L}+I & 0 & 0 \\
			0 & I & 0 \\
			0 & 0 & I
		\end{matrix}
		\right),
	\end{equation*}
 $$G(\xi)=\left(\begin{array}{c}
		-c_{4}\mathcal{L}\bar{x}-\bar{v}\\[0.7em]
		c_{1} g(\bar{x})+c_{2}\mathcal{L}\bar{x}+c_{3} \bar{\lambda}+c_{4} \mathcal{L} \bar{v}+c_{5}\bar{v} \\[0.7em]
		-\mathcal{L} \bar{x}
	\end{array}\right).$$
Then, the algorithm (\ref{eq20}) can be written as  
\begin{equation}\label{eq21}
    \dot{\xi}=-A\xi+F(t).
\end{equation}
where $g(\bar{x})=\nabla f(x)-\nabla f(x^{*})$, $ F(t)=T(\xi(t))=\xi-G(\xi)$. Obviously, $A$ is a positive-definite matrix. 

In the next, we will divide into two steps to prove the convergence of Algorithm \ref{Alg2}. 

(i) In this step, we show that under the assumption of theorem \ref{th2}, $T$ is a contraction, i.e., $\|T(\xi)-T(\xi^{\prime})\|\leq\sigma\|\xi-\xi^{\prime}\|$, where $\sigma=(c_{1}l+c_{3}^{2}+c_{5}^{2}+\sqrt{\left(1+c_{2}^{2}+2c_{4}^{2}\right)\lambda_{\max}\left(L^{\mathrm{T}}L\right)})^{\frac{1}{2}}$.
\begin{align}
\|T(\xi)-T(\xi^{\prime})\|^{2} &= \| \xi-G(\xi)-\xi^{\prime}+G(\xi^{\prime}) \|^{2}\nonumber\\
&=\| \xi-\xi^{\prime}\|^{2}+\|G(\xi)-G(\xi^{\prime})\|^{2} \nonumber\\
&-2\left\langle \xi-\xi^{\prime},G(\xi)-G(\xi^{\prime}) \right\rangle.\label{eq24}
\end{align}
Let $$Q_{1}=\left(\begin{array}{ccc}
	-c_{4}\mathcal{L} & 0 & 0 \\[0.5em]
	c_{2} \mathcal{L} & c_{4}\mathcal{L} & 0 \\[0.5em]
	- \mathcal{L} & 0 & 0
\end{array}\right),Q_{2}=\left(\begin{array}{ccc}
	0 & -I & 0 \\[0.5em]
	0 & c_{5}I & c_{3}I \\[0.5em]
	0 & 0 & 0
\end{array}\right).$$
Consider the second term of the last formula
\begin{align}
    &\left\|G(\xi)-G(\xi^{\prime})\right\|\nonumber \\=& \left\| \begin{array}{c}
		-c_{4}\mathcal{L}\Delta x-\Delta v\\[0.7em]
		c_{1} \Delta \nabla f(x)+c_{2}\mathcal{L}\Delta x+c_{3} \Delta \lambda+c_{4} \mathcal{L} \Delta v+c_{5}\Delta v \\[0.7em]
		-\mathcal{L} \Delta x
	\end{array}\right\|\nonumber\\
    \leq&\left\|\begin{array}{c}
        0\\c_{1}\Delta f\\0
        \end{array}\right\|+ (\left\|Q_{1}\right\|+ \left\|Q_{2}\right\|)\left\|\xi-\xi^{\prime}\right\|\nonumber\\
        \leq&(c_{1}l+\left\|Q_{1}\right\|+ \left\|Q_{2}\right\|)\left\|\xi-\xi^{\prime}\right\|,\label{eq25}
\end{align}
where $\Delta x=\bar{x}-\bar{x}^{\prime}$, $\Delta v=\bar{v}-\bar{v}^{\prime}$, $\Delta \lambda=\bar{\lambda}-\bar{\lambda}^{\prime}$, $\Delta f =\nabla f(\bar{x})-\nabla f(\bar{x}^{\prime})$. The last inequality holds due to the $l$-Lipschitz continuous property of $\nabla f$. From the definition of Kronecker product, we have
\begin{align*}
    \left\|Q_{1}\right\|^{2}&=Q_{1}^{T}Q_{1}=\mathcal{L}^{\mathrm{T}}\mathcal{L}\otimes\left[\begin{array}{ccc}
	c_{2}^{2}+1+c_{4}^{2} & c_{2}c_{4} & 0 \\[0.4em]
	c_{2}c_{4} & c_{4}^{2} & 0 \\[0.4em]
	0 & 0 & 0
\end{array}\right].
\end{align*}

Since $\mathcal{L}^{\mathrm{T}}\mathcal{L}$ is positive semi-definite and
$$\lambda\left(\left[\begin{array}{ccc}
	c_{2}^{2}+1+c_{4}^{2} & c_{2}c_{4} & 0 \\[0.4em]
	c_{2}c_{4} & c_{4}^{2} & 0 \\[0.4em]
	0 & 0 & 0
\end{array}\right]\right)\leq 1+c_{2}^{2}+2c_{4}^{2},$$
it can be concluded that
\begin{equation}\label{eq26}
    \left\|Q_{1}\right\|\leq{\sqrt{\left(1+c_{2}^{2}+2c_{4}^{2}\right)\lambda_{\max}\left(L^{\mathrm{T}}L\right)}}.
\end{equation}
Note that $\left\|Q_{2}\right\|\leq(1+c_{3}^{2}+c_{5}^{2})$. Substituting (\ref{eq26}) into (\ref{eq25}) yields
\begin{align}\label{eq27}
    \left\|G(\xi)-G(\xi^{\prime})\right\|
    \leq&\sigma_{1}\left\|\xi-\xi^{\prime}\right\|,
\end{align}
where $\sigma_{1}=c_{1}l+(1+c_{3}^{2}+c_{5}^{2})+\sqrt{(1+c_{2}^{2}+2c_{4}^{2})\lambda_{\max}(L^{\mathrm{T}}L)}$. From definition \ref{def2}, (\ref{eq27}) can be transformed into
\begin{align}\label{eq28}
    \left\langle G(\xi)-G( \xi^{\prime}), \xi-\xi^{\prime}\right\rangle\geq
    \frac{1}{\sigma_{1}}
    \left\|G(\xi)-G(\xi^{\prime})\right\|^{2}.
\end{align}
Substituting (\ref{eq27}) and (\ref{eq28}) into (\ref{eq24}), we can get
$$\|T(\xi)-T(\xi^{\prime})\|^{2}\leq(\sigma_{1}-1)^{2}\|\xi-\xi^{\prime}\|^{2},$$which implies that $\|T(\xi)-T(\xi^{\prime})\|\leq\sigma\|\xi-\xi^{\prime}\|$.

(ii) In this step, we show the exponentially convergence utilizing the Lyapunov analysis method.

Since $A$ is a positive-definite matrix, all the eigenvalues of $-A$ are located at the open left half plane. Then there exist positive constants $\gamma$ and $\eta$ such that
\begin{equation}\label{eq101}
    \|e^{-A(t-t_{0})}\|\leq\gamma e^{-\eta(t-t_{0})},t\geqslant t_{0}.
\end{equation}
Utilizing the variation of parameter formula, (\ref{eq21}) implies 
\begin{equation*}
    \xi(t)=e^{-A(t-t_{0})}\xi(t_{0})+\int_{t_{0}}^{t}e^{-A(t-s)}F(s)d s,
\end{equation*}
and then
\begin{equation}\label{eq102}
    \|\xi(t)\|\leq \|e^{-A(t-t_{0})}\|\cdot\|\xi(t_{0})\|+\int_{t_{0}}^{t}\|e^{-A(t-s)}\|\|F(s)\|d s.
\end{equation}
By (\ref{eq101}), it follows that
\begin{align}\label{eq103}
    \|\xi(t)\|&\leq\gamma e^{-\eta(t-t_{0}))}\|\xi(t_{0})\|+\gamma\int_{t_{0}}^{t}e^{-\eta(t-s)}\|F(s)\|d s\nonumber\\
    &\leq a e^{-\eta(t-t_{0}))}+b\int_{t_{0}}^{t}e^{-\eta(t-s)}\|\xi(s)\|d s,
\end{align}
where $a=\gamma\|\xi(t_{0})\|$, $b=\gamma\sigma$. Multiply both sides of equation (\ref{eq103}) by $e^{\eta(t-t_{0}))}$, combined with the Gronwall inequality, we have $e^{\eta(t-t_{0}))}\|\xi(t)\|\leq a e^{b(t-t_{0}))}$,
which implies 
\begin{equation}
    \|\xi(t)\|\leq a e^{-(\eta-b)(t-t_{0}))}.
\end{equation}
For $\eta-b\geq 0$, i.e., $\sigma<\frac{\eta}{\gamma}$, $\xi(t)$ can exponentially converge to 0. Then $x_{i}\to{x^{*}}$ as $t\to\infty$, and the convergence rate is estimated by $\eta-b$.

\bibliographystyle{IEEEtran}
\bibliography{pid.bib}

\begin{thebibliography}{10}
\providecommand{\url}[1]{#1}
\csname url@samestyle\endcsname
\providecommand{\newblock}{\relax}
\providecommand{\bibinfo}[2]{#2}
\providecommand{\BIBentrySTDinterwordspacing}{\spaceskip=0pt\relax}
\providecommand{\BIBentryALTinterwordstretchfactor}{4}
\providecommand{\BIBentryALTinterwordspacing}{\spaceskip=\fontdimen2\font plus
\BIBentryALTinterwordstretchfactor\fontdimen3\font minus \fontdimen4\font\relax}
\providecommand{\BIBforeignlanguage}[2]{{%
\expandafter\ifx\csname l@#1\endcsname\relax
\typeout{** WARNING: IEEEtran.bst: No hyphenation pattern has been}%
\typeout{** loaded for the language `#1'. Using the pattern for}%
\typeout{** the default language instead.}%
\else
\language=\csname l@#1\endcsname
\fi
#2}}
\providecommand{\BIBdecl}{\relax}
\BIBdecl

\bibitem{nedic2018network}
A.~Nedi{\'c}, A.~Olshevsky, and M.~G. Rabbat, ``Network topology and communication-computation tradeoffs in decentralized optimization,'' \emph{Proceedings of the IEEE}, vol. 106, no.~5, pp. 953--976, 2018.

\bibitem{yang2019survey}
T.~Yang, X.~Yi, J.~Wu, Y.~Yuan, D.~Wu, Z.~Meng, Y.~Hong, H.~Wang, Z.~Lin, and K.~H. Johansson, ``A survey of distributed optimization,'' \emph{Annual Reviews in Control}, vol.~47, pp. 278--305, 2019.

\bibitem{boyd2011distributed}
S.~Boyd, N.~Parikh, E.~Chu, B.~Peleato, J.~Eckstein \emph{et~al.}, ``Distributed optimization and statistical learning via the alternating direction method of multipliers,'' \emph{Foundations and Trends{\textregistered} in Machine learning}, vol.~3, no.~1, pp. 1--122, 2011.

\bibitem{wang2011control}
J.~Wang and N.~Elia, ``A control perspective for centralized and distributed convex optimization,'' in \emph{2011 50th IEEE conference on decision and control and European control conference}.\hskip 1em plus 0.5em minus 0.4em\relax IEEE, 2011, pp. 3800--3805.

\bibitem{yang2016multi}
S.~Yang, Q.~Liu, and J.~Wang, ``A multi-agent system with a proportional-integral protocol for distributed constrained optimization,'' \emph{IEEE Transactions on Automatic Control}, vol.~62, no.~7, pp. 3461--3467, 2016.

\bibitem{zhang2014distributed}
Y.~Zhang and Y.~Hong, ``Distributed optimization design for second-order multi-agent systems,'' in \emph{Proceedings of the 33rd Chinese control conference}.\hskip 1em plus 0.5em minus 0.4em\relax IEEE, 2014, pp. 1755--1760.

\bibitem{kia2015distributed}
S.~S. Kia, J.~Cort{\'e}s, and S.~Mart{\'\i}nez, ``Distributed convex optimization via continuous-time coordination algorithms with discrete-time communication,'' \emph{Automatica}, vol.~55, pp. 254--264, 2015.

\bibitem{yi2018distributed}
X.~Yi, L.~Yao, T.~Yang, J.~George, and K.~H. Johansson, ``Distributed optimization for second-order multi-agent systems with dynamic event-triggered communication,'' in \emph{2018 IEEE Conference on Decision and Control (CDC)}.\hskip 1em plus 0.5em minus 0.4em\relax IEEE, 2018, pp. 3397--3402.

\bibitem{li2020input}
M.~Li, G.~Chesi, and Y.~Hong, ``Input-feedforward-passivity-based distributed optimization over jointly connected balanced digraphs,'' \emph{IEEE Transactions on Automatic Control}, vol.~66, no.~9, pp. 4117--4131, 2020.

\bibitem{6578120}
B.~Gharesifard and J.~Cortés, ``Distributed continuous-time convex optimization on weight-balanced digraphs,'' \emph{IEEE Transactions on Automatic Control}, vol.~59, no.~3, pp. 781--786, 2014.

\bibitem{yue2021distributed}
D.~Yue, S.~Baldi, J.~Cao, and B.~De~Schutter, ``Distributed adaptive optimization with weight-balancing,'' \emph{IEEE Transactions on Automatic Control}, vol.~67, no.~4, pp. 2068--2075, 2021.

\bibitem{6129483}
J.~Lu and C.~Y. Tang, ``Zero-gradient-sum algorithms for distributed convex optimization: The continuous-time case,'' \emph{IEEE Transactions on Automatic Control}, vol.~57, no.~9, pp. 2348--2354, 2012.

\bibitem{yu2021generalized}
H.~Yu, M.~H. Dhullipalla, and T.~Chen, ``A generalized algorithm for continuous-time distributed optimization,'' in \emph{2021 American Control Conference (ACC)}.\hskip 1em plus 0.5em minus 0.4em\relax IEEE, 2021, pp. 820--825.

\bibitem{borase2021review}
R.~P. Borase, D.~Maghade, S.~Sondkar, and S.~Pawar, ``A review of pid control, tuning methods and applications,'' \emph{International Journal of Dynamics and Control}, vol.~9, pp. 818--827, 2021.

\bibitem{ang2005pid}
K.~H. Ang, G.~Chong, and Y.~Li, ``Pid control system analysis, design, and technology,'' \emph{IEEE transactions on control systems technology}, vol.~13, no.~4, pp. 559--576, 2005.

\bibitem{zhu2022distributed}
W.~Zhu and H.~Tian, ``Distributed convex optimization via proportional-integral-differential algorithm,'' \emph{Measurement and Control}, vol.~55, no. 1-2, pp. 13--20, 2022.

\bibitem{yu2016gradient}
W.~Yu, P.~Yi, and Y.~Hong, ``A gradient-based dissipative continuous-time algorithm for distributed optimization,'' in \emph{2016 35th Chinese Control Conference (CCC)}.\hskip 1em plus 0.5em minus 0.4em\relax IEEE, 2016, pp. 7908--7912.

\bibitem{attouch2000heavy}
H.~Attouch and F.~Alvarez, ``The heavy ball with friction dynamical system for convex constrained minimization problems,'' in \emph{Optimization: Proceedings of the 9th Belgian-French-German Conference on Optimization Namur, September 7--11, 1998}.\hskip 1em plus 0.5em minus 0.4em\relax Springer, 2000, pp. 25--35.

\bibitem{yi2017formation}
X.~Yi, J.~Wei, D.~V. Dimarogonas, and K.~H. Johansson, ``Formation control for multi-agent systems with connectivity preservation and event-triggered controllers,'' \emph{IFAC-PapersOnLine}, vol.~50, no.~1, pp. 9367--9373, 2017.

\bibitem{abdelwahed2013reduced}
S.~Abdelwahed, A.~Asrari, J.~Crider, R.~Dougal, M.~Faruque, Y.~Fu, J.~Langston, Y.~Lee, H.~Mohammadpour, A.~Ouroua \emph{et~al.}, ``Reduced order modeling of a shipboard power system,'' in \emph{2013 IEEE Electric Ship Technologies Symposium (ESTS)}.\hskip 1em plus 0.5em minus 0.4em\relax IEEE, 2013, pp. 256--263.

\bibitem{gonzalez2011effect}
F.~Gonz{\'a}lez, M.~{\'A}. Naya, A.~Luaces, and M.~Gonz{\'a}lez, ``On the effect of multirate co-simulation techniques in the efficiency and accuracy of multibody system dynamics,'' \emph{Multibody System Dynamics}, vol.~25, pp. 461--483, 2011.

\end{thebibliography}

\end{document}